# ON THE ALMOST SURE CONVERGENCE OF SUMS

LUCA PRATELLI AND PIETRO RIGO

ABSTRACT. Two counterexamples, addressing questions raised in [1] and [3], are provided. Both counterexamples are related to chaoses. Let $F_n = Y_n + Z_n$. It may be that $F_n \xrightarrow{a.s.} 0$, $F_n \xrightarrow{L_{2+\delta}} 0$ and $E\{\sup_n |F_n|^\delta\} < \infty$, where $\delta > 0$ and $Y_n$ and $Z_n$ belong to chaoses of uniformly bounded degree, and yet $Y_n$ fails to converge to 0 a.s.

## 1. INTRODUCTION

Throughout, $(X_n : n \in \mathbb{N})$ is a sequence of real independent random variables such that $E(X_n) = 0$ and $E(X_n^2) = 1$. For $p \in \mathbb{N}$, we denote by $\mathcal{F}_p$ the collection of those functions $f : \mathbb{N}^p \to \mathbb{R}$ such that $f$ is symmetric (i.e., invariant under permutations of its arguments), $f$ vanishes on the diagonal (i.e., $f(j_1, \ldots, j_p) = 0$ whenever $j_r = j_s$ for some $r \neq s$), and
$$\sum_{j_1, \ldots, j_p \in \mathbb{N}} f(j_1, \ldots, j_p)^2 < \infty.$$
Moreover, for $p \in \mathbb{N}$ and $f \in \mathcal{F}_p$, we let
$$Q_p(f) = \sum_{j_1, \ldots, j_p \in \mathbb{N}} f(j_1, \ldots, j_p) X_{j_1} \ldots X_{j_p}.$$
Here, the series converges (a.s. and in $L_2$) because of the martingale convergence theorem. In fact,
$$M_n = \sum_{j_1, \ldots, j_p = 1}^{n} f(j_1, \ldots, j_p) X_{j_1} \ldots X_{j_p}$$
is a martingale with respect to the filtration $\sigma(X_1, \ldots, X_n)$ and
$$\sup_n E\{M_n^2\} = \sum_{j_1, \ldots, j_p \in \mathbb{N}} f(j_1, \ldots, j_p)^2.$$
Therefore, $Q_p(f) = \lim_n M_n$ where the equality is a.s. and in the $L_2$-sense.

Multilinear forms in independent random variables, such as $Q_p(f)$, play a role in various frameworks, including multiple stochastic integration, harmonic analysis, Boolean functions, geometry of Banach spaces, random graphs, concentration of measures, Malliavin calculus and the Malliavin-Stein method; see [1] and references therein.

Recently, the following problem has come to the fore. Suppose that
$$F_n \xrightarrow{a.s.} F,$$

---

where $F_n$ and $F$ are real random variables of the form

$$F_n = a_n + \sum_{p=1}^{q} Q_p(f_{n,p}) \quad \text{and} \quad F = a + \sum_{p=1}^{q} Q_p(f_p)$$

for some $q \in \mathbb{N}$, $f_p$, $f_{n,p} \in \mathcal{F}_p$ and $a, a_n \in \mathbb{R}$. Can one deduce that

(1) $\quad a = \lim_n a_n \quad \text{and} \quad Q_p(f_p) \stackrel{a.s.}{=} \lim_n Q_p(f_{n,p}) \quad \text{for all } p = 1, \ldots, q$ ?

As possibly expected, the answer is no. To give conditions for (1), hence, is a quite natural problem.

For instance, by a result of Poly and Zheng [3, Theorem 1.3], condition (1) holds whenever

$$\sup_n E\{|X_n|^{2+\delta}\} < \infty \quad \text{for some } \delta > 0.$$

Subsequently, this result has been improved by [1, Cor. 2.9], where condition (1) is shown to be true whenever $(X_n^2 : n \in \mathbb{N})$ is uniformly integrable.

Next, for our purposes, we also need to recall the following result.

**Theorem 1. (Theorem 1.1 of [3]).** *Let $F_n$ and $F$ be real random variables and $S_q = \bigoplus_{p=0}^{q} \mathcal{H}_p$, where $q \in \mathbb{N}$ and $\mathcal{H}_p$ is the $p$-th Gaussian Wiener chaos (associated to a given isonormal Gaussian process). If $F_n \xrightarrow{a.s.} F$ and $F_n \in S_q$ for all $n$, then*

$$F \in S_q \quad \text{and} \quad J_p(F_n) \xrightarrow{a.s.} J_p(F) \quad \text{for } p = 0, 1, \ldots, q$$

*where $J_p$ is the projection operator onto $\mathcal{H}_p$.*

In connection with the above results, the following questions are raised in [3].

(i) Does condition (1) hold if each $X_n$ is a two point variable ?

(ii) Let $T_q = \bigoplus_{p=0}^{q} \mathcal{I}_p$, where $\mathcal{I}_p$ is the $p$-th Poisson-Wiener chaos; see point (*) below and [2, Chap. 18]. Is it possible to replace $S_q$ with $T_q$ in Theorem 1 ? Precisely, if $F_n \xrightarrow{a.s.} F$ and $F, F_n \in T_q$ for all $n$, can one deduce that $J_p(F_n) \xrightarrow{a.s.} J_p(F)$ for all $p = 0, 1, \ldots, q$ ?

Both (i) and (ii) have been answered (in the negative) by Adamczak; see Examples 2.11 and 3.1 of [1].

It is worth noting that, under the assumptions of Theorem 1, one also obtains $F_n \xrightarrow{L_2} F$ (thanks to the hypercontractivity of the Ornstein-Uhlenbeck semigroup). In addition, by Theorem 3.2 of [1], question (ii) has a positive answer provided the sequence $(F_n)$ is $L_2$-bounded and

(2) $\quad E\{\sup_n |F_n|\} < \infty.$

Whether or not condition (2) can be dropped, however, is an open problem; see [1, Example 3.3]. Hence, questions (i)-(ii) could be restated as follows.

(j) Does condition (1) hold if each $X_n$ is a two point variable and $F_n \xrightarrow{L_2} F$ ?

(jj) Is it possible to remove condition (2) if $F_n \xrightarrow{L_2} F$ ? Precisely, if $F_n \xrightarrow{a.s.} F$, $F_n \xrightarrow{L_2} F$ and $F, F_n \in T_q$ for all $n$, can one deduce that $J_p(F_n) \xrightarrow{a.s.} J_p(F)$ for all $p = 0, 1, \ldots, q$ ?

The only purpose of this note is to answer (in the negative) questions (j)-(jj). Our counterexamples are actually simple (the main tool is the Borel-Cantelli lemma) and help to unify the current state of the art on the subject. The second counterexample suggests that, not only condition (2) can not be removed, but to weaken it seems to be quite hard.

## 2. THE COUNTEREXAMPLES

**Example 2. (Question (j) has a negative answer).** Let $(Y_n)$ be a sequence of independent random variables such that $Y_n \in \{-1, 1\}$ and $p_n = P(Y_n = 1) \in (0, 1)$. Define

$$X_n = \frac{Y_n - 2p_n + 1}{2\sqrt{p_n(1-p_n)}} = \sqrt{\frac{1-p_n}{p_n}} 1_{\{Y_n=1\}} - \sqrt{\frac{p_n}{1-p_n}} 1_{\{Y_n=-1\}},$$

$$F = 0 \quad \text{and} \quad F_n = p_{2n+1} X_{2n} + \sqrt{p_{2n+1}(1-p_{2n+1})} X_{2n} X_{2n+1}.$$

We have to choose $p_n$ in such a way that $F_n \xrightarrow{L_2} 0$, $F_n \xrightarrow{a.s.} 0$, but $p_{2n+1} X_{2n}$ fails to converge to 0 a.s. Take, for instance,

$$p_{2n} = 1/n \quad \text{and} \quad p_{2n+1} = n^{-1/\sqrt{\log n}} \quad \text{for } n \geq 2.$$

Then:

- $E(F_n^2) = p_{2n+1} \to 0$.
- Let $A = \{Y_{2n} = Y_{2n+1} = 1 \text{ for infinitely many } n\}$. Since

$$\sum_n P(Y_{2n} = Y_{2n+1} = 1) = \sum_n n^{-1-1/\sqrt{\log n}} < \infty,$$

then $P(A) = 0$. Moreover, on the set $A^c$, one obtains

$$F_n = F_n 1_{\{Y_{2n+1}=1\}} = F_n 1_{\{Y_{2n}=-1, Y_{2n+1}=1\}} = -\sqrt{1/(n-1)}$$

for each $n$ large enough. (The first equality is because $F_n \neq 0$ if and only if $Y_{2n+1} = 1$). Hence, $F_n \xrightarrow{a.s.} 0$.

- Let $B = \{Y_{2n} = 1 \text{ for infinitely many } n\}$. The Borel-Cantelli lemma yields $P(B) = 1$. Moreover, on the set $B$, one obtains

$$p_{2n+1} X_{2n} = p_{2n+1} \sqrt{\frac{1-p_{2n}}{p_{2n}}} = \sqrt{n-1}\, n^{-1/\sqrt{\log n}}$$

for infinitely many $n$. Therefore, $p_{2n+1} X_{2n}$ fails to converge to 0 a.s. (since $\sqrt{n-1}\, n^{-1/\sqrt{\log n}} \to \infty$).

Example 2 is straightforward. The next example is slightly more elaborate.

**Example 3. (Question (jj) has a negative answer).** Let $(Y_n)$ be a sequence of independent Poisson random variables. Define $\lambda_n = E(Y_n)$ and

$$X_n = \frac{Y_n - \lambda_n}{\sqrt{\lambda_n}}, \quad F = 0, \quad F_n = \lambda_{2n+1} X_{2n} + \sqrt{\lambda_{2n+1}} X_{2n} X_{2n+1}.$$

As shown below, $F_n \in T_2$ and $J_1(F_n) = \lambda_{2n+1} X_{2n}$ for all $n$. Therefore, to answer (jj) in the negative, it suffices to choose $\lambda_n$ in such a way that

$$F_n \xrightarrow{L_2} 0, \quad F_n \xrightarrow{a.s.} 0, \quad \lambda_{2n+1} X_{2n} \text{ fails to converge to 0 a.s.}$$

We will prove a little bit more. In fact, we will also obtain

$$F_n \xrightarrow{L_{2+\delta}} 0 \quad \text{and} \quad E\big\{\sup_n |F_n|^\delta\big\} < \infty \quad \text{for some } \delta > 0.$$

To begin with, we collect some useful facts in a lemma.

**Lemma 4.** *If $Y$ has a Poisson distribution with parameter $\lambda > 0$, then*

$$P(Y > j) \leq \lambda^{j+1}/(j+1)! \quad \text{for } j = 0, 1, \ldots$$

*In particular, $P(Y > 0) \leq \lambda$. Moreover, if $\lambda \in (0, 1]$, then*

$$E\big\{|Y - \lambda|^{5/2}\big\} \leq \sqrt{8}\,\lambda \quad \text{and} \quad E(Y^{5/2}) \leq \sqrt{15}\,\lambda.$$

*Proof.* By Taylor formula, there is $\gamma \in (0, \lambda)$ such that

$$P(Y > j) = e^{-\lambda} \sum_{k > j} \lambda^k/k! = e^{-\lambda}\Big\{e^\lambda - \sum_{k=0}^{j} \lambda^k/k!\Big\} = e^{-\lambda} e^\gamma \lambda^{j+1}/(j+1)!$$

Hence, $P(Y > j) \leq \lambda^{j+1}/(j+1)!$ follows from $\gamma < \lambda$. Next, suppose $\lambda \in (0, 1]$ and note that $E\{|Y - \lambda|\} \leq 2\lambda$ and $E\{(Y - \lambda)^4\} = 3\lambda^2 + \lambda \leq 4\lambda$. Therefore,

$$E\big\{|Y - \lambda|^{5/2}\big\} = E\big\{|Y - \lambda|^2 |Y - \lambda|^{1/2}\big\} \leq \sqrt{E\{(Y - \lambda)^4\} E\{|Y - \lambda|\}} \leq \sqrt{8}\,\lambda.$$

Similarly, since $E(Y^4) = \lambda^4 + 6\lambda^3 + 7\lambda^2 + \lambda \leq 15\lambda$,

$$E(Y^{5/2}) = E\big\{Y^2 Y^{1/2}\big\} \leq \sqrt{E(Y^4)\,E(Y)} \leq \sqrt{15}\,\lambda.$$

□

Next, we note that $F_n = X_{2n} Y_{2n+1}$ and we let

$$\lambda_{2n} = n^{-3/4} \quad \text{and} \quad \lambda_{2n+1} = n^{-5/16}.$$

- First, by Lemma 4, one obtains $F_n \xrightarrow{L_{5/2}} 0$. In fact,

$$E(|F_n|^{5/2}) = E\big\{|X_{2n} Y_{2n+1}|^{5/2}\big\} = \frac{E\big\{|Y_{2n} - \lambda_{2n}|^{5/2}\big\}}{\lambda_{2n}^{5/4}} E(Y_{2n+1}^{5/2})$$

$$\leq \sqrt{120}\,\frac{\lambda_{2n+1}}{\lambda_{2n}^{1/4}} = \sqrt{120}\,n^{-1/8} \longrightarrow 0.$$

- We next prove $F_n \xrightarrow{a.s.} 0$. Fix $\epsilon > 0$ and define

$$A = \big\{Y_{2n+1} > \epsilon\,n^{3/8} \text{ for infinitely many } n\big\} \quad \text{and}$$
$$H = \big\{Y_{2n} \neq 0,\ Y_{2n+1} \neq 0 \text{ for infinitely many } n\big\}.$$

On $(A \cup H)^c = A^c \cap H^c$, for each $n$ large enough, one obtains $Y_{2n+1} = 0$ or $Y_{2n} = 0$ and $Y_{2n+1} \leq \epsilon\,n^{3/8}$. In the first case, $F_n = X_{2n} Y_{2n+1} = 0$, while in the second

$$|F_n| = |X_{2n} Y_{2n+1}| = \sqrt{\lambda_{2n}}\,Y_{2n+1} = n^{-3/8} Y_{2n+1} \leq \epsilon.$$

Hence, to prove $F_n \xrightarrow{a.s.} 0$, it suffices to see that $P(A) = P(H) = 0$. In turn, this follows from

$$\epsilon^2 \sum_n P(Y_{2n+1} > \epsilon\, n^{3/8}) \leq \sum_n \frac{E(Y_{2n+1}^2)}{n^{3/4}} = \sum_n \frac{\lambda_{2n+1} + \lambda_{2n+1}^2}{n^{3/4}}$$

$$\leq 2 \sum_n \frac{\lambda_{2n+1}}{n^{3/4}} = 2 \sum_n n^{-17/16} < \infty$$

and

$$\sum_n P(Y_{2n} \neq 0,\, Y_{2n+1} \neq 0) \leq \sum_n \lambda_{2n} \lambda_{2n+1} = \sum_n n^{-17/16} < \infty.$$

- To show that $\lambda_{2n+1} X_{2n}$ does not converge to 0 a.s., we argue as in Example 2. Let $B = \{Y_{2n} = 1 \text{ for infinitely many } n\}$. The Borel-Cantelli lemma yields $P(B) = 1$. On the other hand, on the set $B$, one obtains

$$\lambda_{2n+1} X_{2n} = \lambda_{2n+1} \frac{1 - \lambda_{2n}}{\sqrt{\lambda_{2n}}} = n^{1/16} - n^{-11/16}$$

for infinitely many $n$. Therefore, $\lambda_{2n+1} X_{2n}$ fails to converge to 0 a.s.

- Finally, we prove $E(M^\delta) < \infty$ for each $\delta \in (0, 1/24)$, where

$$M = \sup_n |F_n|.$$

Since $E(M^\delta) = \int_0^\infty P(M > t^{1/\delta})\, dt$, it suffices to show that

(3) $$P(M > t) = \mathrm{O}(t^{-1/24}) \qquad \text{as } t \to \infty.$$

For each $t > 0$, Lemma 4 yields

$$P(M > t) \leq \sum_n P(|F_n| > t) = \sum_n P(|X_{2n} Y_{2n+1}| > t)$$

$$\leq \sum_n \left\{ P(0 < Y_{2n+1} \leq \sqrt{t},\, |X_{2n}| > \sqrt{t}) + P(Y_{2n+1} > \sqrt{t}) \right\}$$

$$\leq \sum_n P(Y_{2n+1} > 0)\, P(|X_{2n}| > \sqrt{t}) + \sum_n \frac{\lambda_{2n+1}^{[\sqrt{t}]+1}}{([\sqrt{t}]+1)!}$$

Define $a = \sum_n \lambda_{2n+1}^4 = \sum_n n^{-5/4}$. If $t \geq 9$, then $[\sqrt{t}] + 1 \geq 4$ and $\sqrt{\lambda_{2n} t} > \lambda_{2n}$. Hence,

$$P(M > t) \leq \sum_n \lambda_{2n+1} P(|Y_{2n} - \lambda_{2n}| > \sqrt{\lambda_{2n} t}) + \sum_n \frac{\lambda_{2n+1}^4}{[\sqrt{t}]^2}$$

$$\leq \sum_n \lambda_{2n+1} P(Y_{2n} > \sqrt{\lambda_{2n} t}) + \frac{a}{(\sqrt{t}-1)^2}.$$

Let $b = \sum_n \lambda_{2n+1} \lambda_{2n} = \sum_n n^{-17/16}$. Since $\lambda_{2n} = n^{-3/4} \geq t^{-1/2}$ for each $n \leq [t^{2/3}]$, one obtains

$$\sum_n \lambda_{2n+1} P\big(Y_{2n} > \sqrt{\lambda_{2n} t}\big) \leq \sum_{n=1}^{[t^{2/3}]} \lambda_{2n+1} P\big(Y_{2n} > t^{1/4}\big) + \sum_{n>[t^{2/3}]} \lambda_{2n+1} \lambda_{2n}$$

$$\leq \sum_{n=1}^{[t^{2/3}]} \lambda_{2n+1} \lambda_{2n} t^{-1/4} + \sum_{n>[t^{2/3}]} n^{-17/16}$$

$$\leq b t^{-1/4} + 16 [t^{2/3}]^{-1/16} \leq b t^{-1/4} + 16 \{t^{2/3} - 1\}^{-1/16}.$$

To summarize,

$$P(M > t) \leq b t^{-1/4} + 16 \{t^{2/3} - 1\}^{-1/16} + \frac{a}{(\sqrt{t} - 1)^2} \qquad \text{for each } t \geq 9,$$

and this proves condition (3).

To close the paper, it remains to see that $F_n \in T_2$ and $J_1(F_n) = \lambda_{2n+1} X_{2n}$.

**(*) The Poisson-Wiener chaos.** Let $(\mathcal{X}, \mathcal{E}, \mu)$ be a $\sigma$-finite measure space and $N$ a Poisson process with intensity $\mu$ on some probability space $(\Omega, \mathcal{A}, P)$. Then,

$$L_2(\Omega, \sigma(N), P) = \bigoplus_{p=0}^{\infty} \mathcal{I}_p \qquad \text{with } \mathcal{I}_p \perp \mathcal{I}_q \text{ for } p \neq q.$$

The closed subspace $\mathcal{I}_p$, called the *p-th Poisson-Wiener chaos*, is the collection of all random variables of the form $I_p(f)$, where $I_p$ is the multiple Poisson integral of order $p$ and $f$ ranges over the symmetric elements of $L_2(\mathcal{X}^p, \mathcal{E}^p, \mu^p)$. We refer to [2, Chapter 12] for the general definition of $I_p$. Here, it suffices to note that, if $A = A_1 \times \ldots \times A_p$ with $A_i \in \mathcal{E}$, $\mu(A_i) < \infty$ and $A_i \cap A_j = \emptyset$ for $i \neq j$, then

$$I_p\Big(\sum_\pi 1_A \circ \pi\Big) = p! \prod_{i=1}^p \{N(A_i) - \mu(A_i)\}$$

where the sum is over all permutations $\pi$ of $\mathcal{X}^p$.

Next, let $\mathcal{X} = [0, \infty)$, $\mathcal{E}$ the Borel $\sigma$-field, $\mu$ the Lebesgue measure, and

$$Y_n = N(A_n) \qquad \text{where } A_n = \Big(\sum_{i=1}^{n-1} \lambda_i, \sum_{i=1}^n \lambda_i\Big].$$

Then, $(Y_n)$ is an independent sequence, each $Y_n$ has a Poisson distribution with parameter $\lambda_n$, and

$$\sqrt{\lambda_{2n+1}} X_{2n} X_{2n+1} = I_2(f_n) \in \mathcal{I}_2 \qquad \text{where}$$

$$f_n(x,y) = \frac{1}{2\sqrt{\lambda_{2n}}} \big\{1_{A_{2n}}(x) 1_{A_{2n+1}}(y) + 1_{A_{2n}}(y) 1_{A_{2n+1}}(x)\big\}.$$

Similarly,

$$\lambda_{2n+1} X_{2n} = I_1\Big(\lambda_{2n}^{-1/2} \lambda_{2n+1} 1_{A_{2n}}\Big) \in \mathcal{I}_1.$$

Hence,

$$F_n = \lambda_{2n+1} X_{2n} + \sqrt{\lambda_{2n+1}} X_{2n} X_{2n+1} \in T_2 \quad \text{and} \quad J_1(F_n) = \lambda_{2n+1} X_{2n}.$$

Luca Pratelli, Accademia Navale, viale Italia 72, 57100 Livorno, Italy
*E-mail address*: `pratel@mail.dm.unipi.it`

Pietro Rigo (corresponding author), Dipartimento di Scienze Statistiche "P. Fortunati", Università di Bologna, via delle Belle Arti 41, 40126 Bologna, Italy
*E-mail address*: `pietro.rigo@unibo.it`